\documentclass[titlepage,11pt]{article}
\usepackage{latexsym}
\usepackage{amsfonts}
\usepackage{amsmath}
\newcommand\blackslug{\hbox{\hskip 1pt \vrule width 4pt height 8pt depth 1.5pt
        \hskip 1pt}}
\newcommand\bbox{\hfill \quad \blackslug \bigbreak}
\def\d{\hbox{-}}
\def\c{\hbox{-}\cdots\hbox{-}}
\def\ll{,\ldots,}

\title{Disjoint paths in unions of tournaments}
\author{Maria Chudnovsky\thanks{Supported by NSF grant DMS-1265803.}\\
Princeton University, Princeton, NJ 08544, USA
\\
\\
Alex Scott\\
Mathematical Institute, University of Oxford, Oxford OX2 6GG, UK
\\
\\
Paul Seymour\thanks{Supported by NSF grant DMS-1265563 and ONR grant N00014-14-1-0084.}\\
Princeton University, Princeton, NJ 08544, USA}

\date{}

\newtheorem{thm}{}[section]

\newcommand{\Proof}{\noindent{\bf Proof.}\ \ }

\begin{document}
\maketitle
\begin{abstract}
Given $k$ pairs of vertices $(s_i,t_i)\;(1\le i\le k)$ of a digraph $G$, how can we test whether there exist vertex-disjoint directed
paths from $s_i$ to $t_i$ for $1\le i\le k$? This is NP-complete in general digraphs, even for $k = 2$~\cite{FHW},
but in~\cite{tournpaths} we proved that for all fixed $k$, there is a polynomial-time algorithm to solve the problem if $G$ is a tournament
(or more generally, a semicomplete digraph).
Here we prove that for all fixed $k$ there 
is a polynomial-time algorithm to solve the problem when $V(G)$ is partitioned into a bounded number of
sets each inducing a semicomplete digraph (and we are given the partition).
\end{abstract}

\section{Introduction}

A {\em linkage} in a digraph $G$ is a family $L = (P_i\;:1\le i\le k)$ of pairwise vertex-disjoint
directed paths of $G$. 
(With a slight abuse of terminology, we call $k$ the {\em cardinality} of $L$, and $P_1\ll P_k$ its {\em members}.)
Let $s_1,t_1\ll s_k,t_k$ be distinct vertices of a digraph $G$. We call $(G,s_1,t_1\ll s_k,t_k)$ a {\em problem instance}.
A linkage $L = (P_i\;:1\le i\le k)$ in $G$ is {\em for} the problem instance 
if $P_i$ is from $s_i$ to $t_i$ for each $i$.
The {\em $k$ vertex-disjoint paths problem} is to determine whether
there is a linkage for a given problem instance.
Fortune, Hopcroft and Wyllie~\cite{FHW} showed that 
this is NP-complete, even for $k = 2$. 
This motivates the study of subclasses of digraphs for which the problem is polynomial-time solvable. 

In this paper, all digraphs are finite, and without loops or parallel edges; thus if $u,v$ are distinct vertices
of a digraph then there do not exist two edges both from $u$ to $v$, although there may be edges
$uv$ and $vu$. Also, by a ``path'' in a digraph we always mean a directed path.
A digraph is a {\em tournament} if for every pair of distinct vertices $u,v$, exactly one of $uv,vu$ is an edge; and a digraph is
{\em semicomplete} if for all distinct $u,v$, at least one of $uv,vu$ is an edge. Bang-Jensen and Thomassen~\cite{bangtho}
showed:
\begin{thm}\label{NPhard}
The $k$ vertex-disjoint paths problem is NP-complete if $k$ is not fixed, even when $G$ is a tournament. 
\end{thm}

In an earlier paper~\cite{tournpaths} we showed:

\begin{thm}\label{tourthm}
For all fixed $k\ge 1$, the $k$ vertex-disjoint paths problem is solvable in polynomial time if $G$ is semicomplete.
\end{thm}

Can this be extended to more general digraphs? One natural question is, what about digraphs with bounded stability number? 
(A set of vertices is {\em stable}
if no edge has both ends in the set, and the {\em stability number} is the size of the largest stable set.) For
the {\em edge-disjoint} directed paths problem, the bounded stability number case is solvable in polynomial time~\cite{FS}.
But for the vertex-disjoint problem, this extension remains out of our reach; indeed, we suspect the problem might be NP-complete
for digraphs with stability number two.

In this paper we do indeed extend \ref{tourthm} to a wider class of digraphs, motivated also by an application in~\cite{banglars}
where the result of this paper is needed. If $G$ is a digraph, a set $C\subseteq V(G)$
is a {\em clique} of $G$ if the subdigraph of $G$ induced on $C$ is semicomplete. Let us say a {\em clique-partition} for a digraph $G$
is a partition $(C_1\ll C_c)$ of $V(G)$ into cliques (we permit the $C_i$'s to be empty).
Our main result is:

\begin{thm}\label{mainthm}
For all fixed $k$ and $c$, there is a polynomial-time algorithm to solve the $k$ vertex-disjoint directed paths problem
in a digraph $G$ that is given with a clique-partition $(C_1\ll C_c)$. Its running time is $O(|V(G)|^t)$ where 
$t$ is about $4(ck)^5$ for $c,k$ large.
\end{thm}

The idea of the algorithm for \ref{mainthm} is a refinement of that for \ref{tourthm}, presented in~\cite{tournpaths}. 
As before, 
we will define an auxiliary digraph $H$ with two special sets of vertices $S_0, T_0$, and prove that there is a 
path in $H$ from $S_0$ to $T_0$ if and only if there is a linkage for $(G,s_1,t_1\ll s_k,t_k)$. 
Thus to solve the problem of \ref{mainthm} it suffices
to construct $H$ in polynomial time.
In the present context, the auxiliary digraph is more complicated than the one in \cite{tournpaths}, because 
it needs extra bells and whistles to accommodate the
clique-partition of $G$.

There are two extensions of \ref{tourthm} proved in \cite{tournpaths}. First, we were able to determine all the minimal 
$k$-tuples $(x_1\ll x_k)$ such that there is a linkage in which the $i$th path has at most $x_i$ vertices, for $1\le i\le k$. 
We have not been able to do the same in the present context. We can determine the minimum integer $x$ such that there is a linkage 
for the problem instance that uses at most $x$ vertices in total, but we cannot control the individual path lengths.

Let $P$ be a path of a digraph $G$, with vertices $v_1\l v_n$ in order. We say $P$ is {\em minimal} if
$j\le i+1$ for every edge $v_iv_j$ of $G$ with $1\le i,j\le n$.
We also showed in \cite{tournpaths} that essentially the same algorithm works for ``$d$-path-dominant'' digraphs instead of just
semicomplete digraphs (these are digraphs in which every $d$-vertex minimal path contains a neighbour of every vertex). Again, we were not able to extend this to
the present context.

\section{The quest for an auxiliary digraph}

Our method is to define an auxiliary digraph $H$, with two special sets of vertices $S_0, T_0$, in such a way that
there is a path in $H$ between $S_0$ and $T_0$ if and only if a linkage exists for $(G, s_1,t_1\ll s_k,t_k)$.
We refer to the parts of this statement as the ``if'' direction and the ``only if'' direction.
To make a polynomial-time algorithm, we need that (a) the number of vertices of $H$
is at most some polynomial in $|V(G)|$, and (b) we can construct $H$ in polynomial time without knowledge of a linkage in $G$.

Here are some attempts, to explain the difficulty and the way we solve it. First, we might try: let $V(H)$ be the set of all $k$-tuples
of distinct vertices of $G$; let $S_0$ contain just the $k$-tuple $(s_1\ll s_k)$, and define $T_0$ similarly; and say vertex
$(u_1\ll u_k)$ of $H$ is adjacent in $H$ to vertex $(v_1\ll v_k)$ if $u_i = v_i$ for all $i$ except one, and
$v_i$ is adjacent from $u_i$
for the exceptional value. We can certainly construct this in polynomial time; and it is easy to see that ``if'' direction holds;
but the ``only if'' direction fails. There might be a path from $S_0$ to $T_0$ in $H$, for which when we trace
out the trajectory in $G$ of the $i$th coordinate, we obtain a walk from $s_i$ to $t_i$ rather than a path (not a problem, we could
short-cut); but worse, the trajectory of one coordinate might use vertices that also have been used by the trajectory of another coordinate.
This is the main difficulty; how can we avoid it?

If $L = (P_i\;:1\le i\le k)$ is a linkage, we define $V(L)$ to be $V(P_1)\cup\cdots\cup V(P_k)$.
A second attempt: let us try to somehow mark the vertices that have already been used, so that they do not get used twice.
Let $H$ consist of $k+1$-tuples, in which the first $k$ terms are vertices of $G$ and the last is a subset of $V(G)$. Say
$(v_1\ll v_k, D)$ is adjacent to $(v_1'\ll v_k', D')$ if again $v_i = v_i'$ for all $i$ except one, and for the last value of $i$,
$v_i'$ is adjacent from $v_i$, and $v_i\notin D$, and $D\cup \{v_i\}= D'$. Take 
$S_0 = \{(s_1\ll s_k, \{s_1\ll s_k\})\}$
and $T_0$ to be all terms of the form $\{(t_1\ll t_k, D)\}$. Then both ``if'' and ``only if'' directions works;
but $H$ has exponential size.

This is of course still naive in several ways. One is that, if the linkage exists, we are tracing it out by
walking $k$-tuples of vertices along its paths,
but not being clever about the sequence of moves of these $k$-tuples. We don't need every sequence of moves of $k$-tuples
that traces out the linkage to correspond to a path in our auxiliary digraph $H$ -- one such sequence giving a path of $H$
would be enough -- so we are being
wasteful here. We could afford to remove some parts of $H$ to make it smaller,
as long as we preserve the property that every linkage in $G$ gives us at least one path in $H$. And even this is wasteful -- 
we don't need every linkage for the problem instance to give a path;
we might as well just make sure that the linkages $L$ work that have vertex set $V(L)$ as small as possible.
These ``minimum'' linkages are nicer than general ones, so this helps.

Suppose we could generate a set $\mathcal{D}$ of polynomially many subsets of $V(G)$, with the following property: that for every 
minimum linkage $L$, there is a way of tracing $L$ with $k$-tuples such that at each stage, the set of vertices that
have been used already is a member of $\mathcal{D}$. This would be ideal, because then we make an auxiliary digraph with
vertices of the form $(v_1\ll v_k,D)$ and adjacency as before, but only using sets $D\in \mathcal{D}$, and
this would all work. However, in general there is no such set $\mathcal{D}$; even for $k = 1$ it is easy to see that there are
tournaments in which every set $\mathcal{D}$ with the property above is exponentially large.

But we are getting closer to an answer. Suppose we could find a polynomially-sized set of subsets $\mathcal{D}$, with the
property that for every minimum linkage $L$, there is a way of tracing out $L$ with $k$-tuples of vertices, such that for every
$k$-tuple $(v_1\ll v_k)$ used in this tracing, there is a set $D\in \mathcal{D}$ which includes the vertices already used,
and includes none of those in the remainder of $V(L)$ (and possibly contains some vertices not used by the linkage). As far as we see,
this would
{\em not} yet be enough, because there seems no way to define the auxiliary graph. We would take
$V(H)$ to be the set of all $k+1$-tuples $(v_1\ll v_k,D)$ where $v_1\ll v_k\in V(G)$ and $D\in \mathcal{D}$, but how should we 
define adjacency in $H$? If $(v_1\ll v_k,D)$ is to be adjacent to $(v_1'\ll v_k',D')$ in $H$, we would presumably want at least that
\begin{itemize}
\item $v_1\ll v_k\in D'$
\item $v_i = v_i'$ for all values of $i$ except one; and
\item $v_i$ is adjacent
to $v_i'$ and $v_i'\in D'\setminus D$ for the
exceptional value of $i$.
\end{itemize}
If we make this the definition of adjacency in $H$ then ``if'' direction works, but the ``only if'' direction fails. If we impose
the additional condition 
\begin{itemize}
\item $D\subseteq D'$
\end{itemize}
then the ``only if'' direction works, but the ``if'' direction fails.

To make the ``if'' direction work (for the four-bullet version of $H$ described above), we need $\mathcal{D}$ to have the following 
additional property:
that, for each $k$-tuple $(v_1\ll v_k)$ used to trace a minimum linkage $L$, there exists $D\in \mathcal{D}$
that intersects $V(L)$ in the set of vertices already used, such that each set $D$ is a subset of the next.
(This used to be automatic when $D$ was just the set of vertices that has been used already; but now that $D$ may contain vertices 
not in $V(L)$, we must impose it as an extra condition.)
That then would work. There {\em is} indeed such a set $\mathcal{D}$
when $G$ is a semicomplete digraph, and that was the idea of our algorithm in \cite{tournpaths}.
Unfortunately, in the present case all we know is that $G$ admits a clique-partition into a bounded number of cliques,
and we have not been able to prove that such a set $\mathcal{D}$ exists, and suspect that in general it does not.

Let us stop trying to trace out the linkage with $k$-tuples of vertices, and trace it out in a different way,
suggested by \ref{vertexorder}. 
That lemma, the key result of the paper, provides, for any minimum linkage $L$, an enumeration of the vertices in $V(L)$, 
which has some useful properties. It gives a sequence of subsets of $V(L)$, starting from $\{s_1\ll s_k\}$ and growing one
vertex at a time until it reaches $V(L)\setminus \{t_1\ll t_k\}$; and each path of the minimum linkage winds in and out of 
any set in this sequence only a bounded number of times. (The enumeration has some other useful properties too that
will be introduced later.) 
We have therefore a sequence of partitions $(A_h,B_h)\; (h =0\ll n) $ of $V(L)$;
and for each $(A_h,B_h)$, there are only at most constantly many (at most $K$ say) ``jumping'' edges
(edges of the linkage paths that pass from $A$ to $B$
or from $B$ to $A$). Let $J_h$ be the set of jumping edges at stage
$h$; then we can regard the sets $J_h\;(0\le h\le n)$ as tracing out the linkage (albeit not as nicely as before: at a general
stage $h$ we will have traced some disjoint set of subpaths of each member of $L$, not just one initial subpath).
Let us
try to design an auxiliary digraph $H$  with the sets $J_h$ replacing the $k$-tuples of vertices. When there is a minimum linkage
$L$ in $G$, and we take sets of jumping edges $J_h\;(0\le h\le n)$ tracing it, the corresponding vertex of $H$ at stage $h$
will be the pair $(J_h,D_h)$.
We need $D_h$ to have three properties:
\begin{itemize}
\item $D_h$ must contain the vertices already used by the partial tracing of the linkage $L$ (that is, $B_h\subseteq D_h$), and
must not contain any vertices in $A_h$ (but it is allowed to contain vertices not in $V(L)$);
\item as $h$ increases, each set $D_h$ must be a subset of the next; and
\item there must be a polynomially-size set $\mathcal{D}$ of subsets of $V(G)$, containing all the sets $D_h$ produced by the chosen
tracing of
the minimum linkage. The sets $D_h$ depend on the choice of $L$; but, crucially,
we must be able to define $\mathcal{D}$ without knowledge of $L$.
\end{itemize}
It will follow from the other desirable features of \ref{vertexorder} that $\mathcal{D}$ and the sets $D_h$ exist with these three properties.
Then we define $H$ to be the digraph with vertex set all the pairs $(J,D)$ where $J$ is a set of at most $K$ edges and
$D\in \mathcal{D}$, and define adjacency in the natural way, and it nearly works;
the problem is, a path in $H$ yields a linkage in $G$ with $k$ paths, all starting in $\{s_1\ll s_k\}$ and ending in $\{t_1\ll t_k\}$,
but not necessarily linking $s_i$ to $t_i$ for $1\le i\le k$. This used not to be a problem because we used to have $k$-tuples
of vertices,
so we could tell which vertex was supposed to belong to which path; but now we are tracing the linkage with sets of edges,
and we can't tell any more which edge is supposed to be in which path. We can fix this by partitioning each set of edges into
$k$ labelled subsets and redefine the adjacency in $H$ to respect the partitions; in other words,
trace with sets of {\em coloured} edges, where the colours are $1\ll k$, and we can tell from the colour of an edge which
path it belongs to.
Doing all this in detail is the content of the remainder of the paper.

\section{The key lemma}

The reduction of the linkage question to the question about finding one path in a different digraph is
thus a more-or-less straightforward consequence of \ref{vertexorder}, and this section is to prove that lemma.
We need a few definitions first.
If $P$ is a directed path of a digraph $G$, its {\em length} is $|E(P)|$ (every path has at least one vertex); 
and $s(P), t(P)$ denote the first and last vertices of $P$, respectively.
If $F$ is a digraph and $v\in V(F)$, $F\setminus v$ denotes the digraph obtained from $F$ by deleting $v$; and
if $X\subseteq V(F)$,
$F[X]$ denotes the subdigraph of $F$ induced on $X$, and $F\setminus X$ denotes the subdigraph obtained by deleting all vertices in $X$.

Now let $L = (P_i\;:1\le i\le k)$ be a linkage in $G$. The linkage
$L$ is {\em minimum} if there is no linkage $L' = (P_i'\;:1\le i\le k)$ in $G$ with $|V(L')|<|V(L)|$ 
joining the same $k$ pairs of vertices (that is, such that $s(P_i) = s(P_i')$ and $t(P_i) = t(P_i')$ for $1\le i\le k$).  
A vertex $v$ is an {\em
internal vertex} of $L$ if $v\in V(L)$, and $v$ is not at either end of any member of $L$. A linkage $L$ is {\em internally disjoint} from a linkage $L'$
if no internal vertex of $L$ belongs to $V(L')$ (note that this does not imply that $L'$ is internally disjoint from $L$); 
and we say that $L,L'$ are {\em internally disjoint} if each of them is internally disjoint from the other (and thus all vertices in $V(L)\cap V(L')$  must be ends of
paths in both $L$ and $L'$)

Let $Q,R$ be paths of a digraph $G$. A {\em planar $(Q,R)$-matching}
is a linkage $(M_j\;:1\le j\le n)$ for some $n\ge 0$ (and we call $n$ its {\em cardinality}), such that 
\begin{itemize}
\item $M_1\ll M_n$ each have at most three vertices;
\item $s(M_1)\ll s(M_n)$ are vertices of $Q$, in order in $Q$; and
\item $t(M_1)\ll t(M_n)$ are vertices of $R$, in order in $R$.
\end{itemize}
(It is convenient not to insist that $Q,R$ are vertex-disjoint; but in all our applications, the planar matching will
be between subpaths $Q', R'$ of $Q,R$ respectively that {\em are} vertex-disjoint.)
If $X,Y\subseteq V(G)$, and each $M_j$ has first vertex in $X$ and last vertex in $Y$, we say this planar
$(Q,R)$-matching is {\em from $X$ to $Y$}.

If $P$ is a directed path, a subpath $Q$ of $P$ with $s(Q)=s(P)$ is called an {\em initial} subpath.
Let $L = (P_1\ll P_k)$ be a linkage for a problem instance $(G,s_1,t_1\ll s_k, t_k)$. 
Let $C\subseteq V(G)$ be a clique.
A subset $B\subseteq C$ is said to be {\em $C$-acceptable}
(for $L$) if (where $A=C\setminus B$):
\begin{itemize}
\item $\{s_1\ll s_k\}\cap C\subseteq B$ and $\{t_1\ll t_k\}\cap B =\emptyset$;
\item for all $i\in \{1\ll k\}$, there is an initial subpath $Q$ of $P_i$ with $V(Q)\cap C=V(P_i)\cap B$; and
\item for all $i,j\in \{1\ll k\}$, there is no planar $(P_i, P_j)$-matching $L'$ from $B$ to $A$ of cardinality $k^2+k+2$ 
internally disjoint from $L$.
\end{itemize}
The next result is a modification of theorem 2.1 of~\cite{tournpaths}.

\begin{thm}\label{C-acceptable}
Let $(G,s_1,t_1\ll s_k,t_k)$ be a problem instance, 
and let $L = (P_1\ll P_k)$ be a minimum linkage for $(G,s_1,t_1\ll s_k,t_k)$.
Let $C$ be a clique of $G$, and
suppose that $B\subseteq V(L)$ is $C$-acceptable for $L$ and $B\ne (V(L)\cap C)\setminus \{t_1\ll t_k\}$. 
Then there exists $v\in (V(L)\cap C)\setminus (B\cup \{t_1\ll t_k\})$ such that 
$B\cup \{v\}$ is $C$-acceptable for $L$.
\end{thm}
\Proof Let $A = C\setminus B$. For $1\le i\le k$, let $r_i$ be the first vertex of $P_i$ in $A\setminus \{t_i\}$, if there
is such a vertex; and let $q_i$ be the vertex immediately preceding it in $P_i$.
Since $L$ is a minimum linkage, we have:
\\
\\
(1) {\em For $1\le i\le k$, $P_i$ is a minimal path of $G$, and in particular, if $r_i$ exists then 
the only edge of $G$ from $V(P_i)\cap B$ to $V(P_i)\cap A$ (if there is one) is $q_ir_i$. Moreover, 
every three-vertex path from $V(P_i)\cap B$
to $V(P_i)\cap A$ with internal vertex in $V(G)\setminus V(L)$ uses at least one of $q_i, r_i$. Consequently,
there is no planar $(P_i, P_i)$-matching from $B$ to $A$ of cardinality three internally disjoint from $L$.}

\bigskip
From (1), the theorem holds if $k = 1$, setting $v=r_1$, so we may assume that $k\ge 2$.
\\
\\
(2) {\em We may assume that for all $i\in \{1\ll k\}$, if $r_i$ exists then
for some $j\in \{1\ll k\}\setminus \{i\}$, $r_j$ exists and there is a $(P_i, P_j)$-planar matching from $B$ to 
$A\setminus \{r_j\}$ 
of cardinality $k^2+k$ internally disjoint from $L$.}
\\
\\
For let $i\in \{1\ll k\}$ such that $r_i$ exists.
We may assume that $B\cup \{r_i\}$ is not $C$-acceptable. Consequently, one of the three conditions in the definition of
``$C$-acceptable'' is not satisfied by $B\cup \{r_i\}$. The first two are satisfied since $r_i$ is the first vertex of $P_i$ in 
$C\setminus B$ and $r_i\ne t_i$.
Thus the third is false, and so for some $i',j\in \{1\ll k\}$, 
there is a planar $(P_{i'}, P_{j})$-matching from $B\cup \{r_i\}$ to $A\setminus \{r_i\}$ of cardinality $k^2+k+2$ 
internally disjoint from $L$.
Since there is no such matching from $B$ to $A$, 
it follows that $i' = i$, and $r_j$ exists, and 
there is a planar $(P_{i}, P_{j})$-matching from $B$ to $A\setminus \{r_j\}$ of cardinality $k^2+k$ internally disjoint from $L$.
Since $k^2+k\ge 4$ (because $k\ge 2$), (1) implies that $j\ne i$. This proves (2).
\\
\\
(3) {\em We may assume that for some $p\ge 2$, and for all $i$ with $1\le i < p$, there is a planar $(P_{i}, P_{i+1})$-matching 
from $B$ to $A\setminus \{r_{i+1}\}$
of cardinality $k^2+k$ internally disjoint from $L$,
and there is a planar $(P_{p}, P_{1})$-matching from $B$ to $A\setminus \{r_1\}$ of cardinality $k^2+k$ internally disjoint from $L$.}
\\
\\
For since $B\ne C\setminus \{t_1\ll t_k\}$, there exists $i\in \{1\ll k\}$ such that $r_i$ exists. 
By repeated application of (2), there exist 
$p\ge 2$ and distinct $h_1\ll h_p\in\{1\ll k\}$ such that for $1\le i\le p$
there is a planar $(P_{h_i}, P_{h_{i+1}})$-matching from $B$ to $A\setminus \{r_{h_{i+1}}\}$ of cardinality $k^2+k$ internally disjoint from $L$, where $h_{p+1} = h_1$.
Without loss
of generality, we may assume that $h_i = i$ for $1\le i\le p$. This proves (3).

\bigskip
Let us say a planar $(Q,R)$-matching is {\em $2$-spaced} if no edge of $Q$ or $R$ meets more than one member of
the matching.
\\
\\
(4) {\em We may assume that for some $p\ge 2$, and for all $i$ with $1\le i < p$, there is a planar $(P_{i}, P_{i+1})$-matching $L_i$ 
from $B$ to $A\setminus \{r_{i+1}\}$,
and there is a planar $(P_{p}, P_{1})$-matching $L_p$ from $B$ to $A\setminus \{r_1\}$, such that 
\begin{itemize}
\item $L_1\ll L_p$ all have cardinality $p$
\item they are pairwise internally disjoint
\item each of $L_1\ll L_p$ is internally disjoint from $L$, and
\item each of $L_1\ll L_p$  is $2$-spaced.
\end{itemize}
}
For let $L_i'$ be a planar $(P_{i}, P_{i+1})$-matching from $B$ to $A\setminus \{r_{i+1}\}$ of cardinality $k^2+k$ 
internally disjoint from $L$, for $1\le i<p$, and
let $L_p'$ be a planar $(P_{p}, P_{1})$-matching from $B$ to $A\setminus \{r_1\}$ of cardinality $k^2+k$ 
internally disjoint from $L$. We choose $L_i\subseteq L_{i}'$ inductively.
Suppose that for some $h<p$, we have chosen $L_1\ll L_h$, such that
\begin{itemize}
\item $L_1\ll L_h$ all have cardinality $p$
\item they are pairwise internally disjoint, and
\item each of $L_1\ll L_h$  is $2$-spaced.
\end{itemize}
We define $L_{h+1}$ as follows. The union of the sets of internal vertices of $L_1\ll L_h$ has cardinality at most $hp\le k(k-1)$, and
so $L_{h+1}'$ includes a planar $(P_{h+1}, P_{h+2})$-matching from $B$ to $A\setminus \{r_{h+2}\}$ (or $(P_{p}, P_{1})$-matching 
from $B$ to $A\setminus \{r_1\}$, if $h = p-1$) 
of cardinality $k^2+k-k(k-1)= 2k$, internally disjoint from each of $L_1\ll L_h$.
By ordering the members of this matching in their natural order, and taking only the $i$th terms where $i<2p$ is odd, 
we obtain a $2$-spaced
matching of cardinality $p$. Let this be $L_{h+1}$. This completes the inductive definition of $L_1\ll L_p$, and so proves (4).

\bigskip
Henceforth we read subscripts modulo $p$.
For $1\le i \le p$, let $L_i = \{M^1_i\ll M^p_i\}$, numbered in order; thus, if $q^h_i$ and $r^h_{i+1}$ denote the first and last vertices of $M^h_i$,
then $q^1_i\ll q^p_i$ are distinct and in order in $Q_i$, and $r_{i+1}, r^1_{i+1}\ll r^p_{i+1}$ are
distinct and in order in $R_{i+1}$.

Since $r^{h}_{i+1}\ne r_{i+1}$, (1) implies that $r^{h}_{i+1}$ is not adjacent from 
$q^{h+1}_{i+1}$; and so there is an edge from $r^{h}_{i+1}$ to $q^{h+1}_{i+1}$, since $C$ is semicomplete.
For $1\le i\le p$, and $1\le h<p$, let $S^h_i$ be the path
$$q^h_i\d M^h_i\d r^h_{i+1}\d q^{h+1}_{i+1};$$
then $S^h_i$ is a path from $q^h_i$ to $q^{h+1}_{i+1}$, of length at most $3$.
Thus concatenating 
$S^1_i, S^2_{i+1}\ll S^{p-1}_{i+p-2}$ and $M^p_{i+p-1}$ gives a path $T_i'$ from $q^1_i$ to $r^p_{i}$ of length at most
$3p-1$ (since $M^p_{i+p-1}$ has at most three vertices, from the definition of a planar matching). 
The subpath $T_i$ of $P_i$ from $q^1_i$ to $r^p_{i}$ has length at least $4(p-1)+2$, since $L_{i-1}, L_i$ are $2$-spaced and 
$r_i$ is different from $r^1_i$; and so $T_i$ has length strictly greater than that of $T_i'$.  Let $P_i'$ be obtained from $P_i$ by replacing the
subpath $T_i$ by $T_i'$, for $1\le i\le p$, and let $P_{i'} = P_i$ for $p+1\le i\le k$.
Then $\{P_1'\ll P_k'\}$ is a linkage for $(G,s_1,t_1\ll s_,t_k)$, contradicting that $L$ is minimum.
This proves \ref{C-acceptable}.~\bbox

Let $P$ be a path of a digraph $G$, and let $X,Y$ be disjoint subsets of $V(G)$. Let $v_1\ll v_t$ be distinct vertices of $P$, in order in $P$.
This sequence is {\em $(X,Y)$-alternating} if $t$ is even and $v_i\in X$ for $i$ odd and $v_i\in Y$ for $i$ even. 
The {\em $(X,Y)$-wiggle number} of $P$
is half the length of the longest $(X,Y)$-alternating sequence $v_1\ll v_t$ where $v_1\ll v_t$ are in order in $P$.
Next we need a lemma, the following:

\begin{thm}\label{disjsegments}
Let $w>0$, let $L$ be a linkage in $G$, and let $Q_1\ll Q_c$ each be a subpath of some member of $L$. Let $X_1,Y_1,X_2,Y_2\ll X_c, Y_c$ be
pairwise disjoint subsets of $V(L)$, such that the $(X_i,Y_i)$-wiggle number of $Q_i$ is at least $cw$ for $1\le i\le c$. 
Then for $1\le i\le c$ there is a subpath $R_i$ of $Q_i$, such that the  $(X_i,Y_i)$-wiggle number of $R_i$ is at least $w$, and
the paths $R_1\ll R_c$ are pairwise vertex-disjoint.
\end{thm}
\Proof We proceed by induction on $c$. If $c=1$ the result holds, so we assume that $c\ge 2$. 
Choose an initial subpath $P_0$ of some member of $L$, minimal such that for some $i\in \{1\ll c\}$, $P_0\cap Q_i$ is nonnull and
the $(X_i,Y_i)$-wiggle number of the path $P_0\cap Q_i$ is at least $w$. We may assume that $i=c$, that is, the
$(X_c,Y_c)$-wiggle number of $P_0\cap Q_c$ is at least $w$. Let $R_c = P_0\cap Q_c$.
For $1\le i<c$ let $Q_i'=Q_i\setminus V(P_0)$. From the minimality of $P_0$, the $(X_i,Y_i)$-wiggle number of $P_0\cap Q_i$
is at most $w$, and either this number is less than $w$ or the last vertex of $P_0$ is in $Y_i$. So 
the $(X_i,Y_i)$-wiggle number of $Q_i'$ is at least $w(c-1)$. The result follows from 
the inductive hypothesis applied
to $Q_1'\ll Q_{c-1}'$, since they are all disjoint from $R_c$.~\bbox

Let $k,c\ge 1$ and define $z = c(c(k^2+k+1)+ k+2)$.
Now let $L=(P_1\ll P_k)$ be a minimum linkage for a problem instance $(G,s_1,t_1\ll s_k,t_k)$, and let $(C_1\ll C_c)$ be a clique-partition of $G$.
Let $B\subseteq V(G)$ and $A=V(G)\setminus B$. We say that $B$ is {\em acceptable} if:
\begin{itemize}
\item $s_1\ll s_k\in B$ and $t_1\ll t_k\notin B$;
\item for $1\le a\le c$, $B\cap C_a$ is $C_a$-acceptable; and
\item for all distinct $a,b$ with $1\le a,b\le c$, and for $1\le i\le k$, the $(B\cap C_b, A\cap C_a)$-wiggle number of $P_i$ is
at most $z$.
\end{itemize}

\begin{thm}\label{acceptable}
Let $k,c,z$ and $(G,s_1,t_1\ll s_k,t_k)$, $L$ and $(C_1\ll C_c)$ be as above. Let $B\subseteq V(L)$
be acceptable, with $B\ne V(L)\setminus \{t_1\ll t_k\}$. 
Then there exists $v\in V(L)\setminus (B\cup \{t_1\ll t_k\})$ such that $B\cup \{v\}$ is acceptable.
\end{thm}
\Proof 
Let $w = c(k^2+k+1)+ k+2$.
Let $A=V(L)\setminus B$. For $1\le a\le c$, if $C_a\cap A\not \subseteq \{t_1\ll t_k\}$, choose $r_a\in (C_a\cap A)\setminus \{t_1\ll t_k\}$ such that 
$B\cap C_a\cup \{r_a\}$ is $C_a$-acceptable (this is possible by \ref{C-acceptable}). 
Since $B\ne V(L)\setminus\{t_1\ll t_k\}$, there is at least one value of 
$a\in \{1\ll c\}$ such that $r_a$ exists. Suppose that there is no $a\in \{1\ll c\}$ 
such that $r_a$ exists and $B\cup \{r_a\}$ is acceptable.
\\
\\
(1) {\em If $1\le a\le c$ and $r_a$ exists, 
let $r_a\in V(P_i)$; then there exists $b\in \{1\ll c\}$ with $b\ne a$
such that the $(B\cap C_a, A\cap C_b)$-wiggle number of $P_i$ is at least $z$.}
\\
\\
Let $B'=B\cup \{r_a\}$. From the choice of $r_a$, it follows that 
$B\cap C_b$ is $C_b$-acceptable for $1\le b\le c$; and so, since $B'$ is not acceptable,
there exist distinct $a',b'\in \{1\ll c\}$, and $i\in \{1\ll k\}$,
such that the $(B'\cap C_{b'}, A'\cap C_{a'})$-wiggle number of $P_i$ is
at least $z+1$, where $A' = A\setminus \{r_a\}$. Let $v_1,v_2\ll v_{2z+2}$ be a $(B'\cap C_{b'},A'\cap C_{a'})$-alternating  
sequence of vertices of $P_i$, in order
in $P_i$. This sequence is not $(B\cap C_{b'},A\cap C_{a'})$-alternating,
since $B$ is acceptable; and so one of $v_1\ll v_{2z+2}$ equals $r_a$. In particular,
$r_a$ belongs to one of $A'\cap C_{a'}$, $B'\cap C_{b'}$. Since $r_a\notin A'$, it follows that $r_a\in B'\cap C_{b'}$,
and so $a=b'$. Since $B$ is $C_{a}$-acceptable, we deduce that $r_a$ is later in $P_i$
than every vertex of $P_i$ in $B\cap C_a$; and since $v_1\ll v_{2z+2}$ are in order in $P_i$, and 
$$v_{2z+1}\in B'\cap C_{b'}= (B\cup \{r_a\})\cap C_a,$$ 
it follows that $r_a=v_{2z+1}$. Consequently $v_1\ll v_{2z}$ is $(B\cap C_{a},A\cap C_{a'})$-alternating,
and so setting $b=a'$ satisfies the claim. This
proves (1).

\bigskip

From (1), 
and \ref{disjsegments}, for each $a\in \{1\ll c\}$ such that $r_a$ exists, there is a subpath $R_a$ of some member of $L$
and $b\in \{1\ll c\}$ with $b\ne a$ such that the $(B\cap C_a, A\cap C_b)$-wiggle number of $R_a$ is at least $w$, and
the paths $R_a\;(1\le a\le c)$ are pairwise disjoint (if they exist). In particular, if $b$ is as above then
$C_b\cap A\not \subseteq \{t_1\ll t_k\}$ and so $r_b$ exists. Renumbering, we may assume that for some $p$ with $2\le p\le c$:

\begin{itemize}
\item there are paths $R_1\ll R_p$, each a subpath of some member of $L$ and pairwise disjoint;
\item for $1\le a < p$, the $(B\cap C_a,A\cap C_{a+1})$-wiggle number of $R_a$ is at least $w$,
and the $(B\cap C_p, A\cap C_{1})$-wiggle number of $R_p$ is at least $w$.
\end{itemize}

Consequently the $(A\cap C_{a+1},B\cap C_a)$-wiggle number of $R_a$ is at least $w-1$.
For $1\le a\le p$, choose vertices $x^1_a,y^1_a\ll x^{w-1}_a,y^{w-1}_a$ in order in $R_a$ and $(A\cap C_{a+1},B\cap C_a)$-alternating 
(henceforth we read subscripts modulo $p$).
Thus $x^1_a, x^2_a\ll x^{w-1}_a\in A\cap C_{a+1}$, and $y^1_a, y^2_a\ll y^{w-1}_a\in B\cap C_a$.
Since $B$ is $C_{a}$-acceptable, there is no planar $(R_{a}, R_{a-1})$-matching of cardinality $k^2+k+2$ from
$B\cap C_{a}$ to $A\cap C_{a}$ internally disjoint from $L$; and in particular, since $x^1_{a-1}, x^2_{a-1}\ll x^{w-1}_{a-1}\in A\cap C_a$
and $y^1_{a}, y^2_{a}\ll y^{w-1}_{a}\in B\cap C_{a}$, it follows that $y^i_{a}$ is adjacent to $x^{i}_{a-1}$ for at most
$k^2+k+1$ values of $i$. Since $C_{a}$ is semicomplete, it follows that $y^i_{a}$ is adjacent from $x^{i}_{a-1}$ for all except
$k^2+k+1$ values of $i$. Hence there exists $I\subseteq \{1\ll w-1\}$ of cardinality at least $w-1-c(k^2+k+1)=k+1$, such that
$y^i_{a}$ is adjacent from $x^{i}_{a-1}$
for all $i\in I$ and $a\in \{1\ll p\}$. Renumbering, and using the fact that $k+1\ge p+1$, 
it follows that for $1\le a\le p$,
there exist $u^1_a,v^1_a\ll u^{p}_a,v^{p}_a$ in order in $R_a$ and $(A\cap C_{a+1}, B\cap C_a)$-alternating,
such that $v^i_{a}$ is adjacent from $u^i_{a-1}$
for all $i,a$ with $1\le i\le p$ and $1\le a\le p$, and in addition $u^1_1, v^1_1$ are not consecutive vertices of $R_1$.

For $1\le a\le p$ and $1\le i< p$, let $T^i_a$ be the subpath of $R_i$ with first vertex $v^i_a$ and last vertex $u^{i+1}_a$.
Then for $1\le a\le p$,
$$u^1_a\d v^1_{a+1}\d T^1_{a+1}\d u^2_{a+1}\d v^2_{a+2}\d T^2_{a+2}\c T^{p-1}_{a-1}\d u^p_{a-1}\d v^p_a$$
is a directed path from $u^1_a$ to $v^p_a$, say $M_a$. For $1\le a\le c$ let $R_a'$ be the subpath of $R_a$ from $u^1_a$ to $v^p_a$.
The paths $M_1\ll M_p$ are pairwise disjoint, and 
$$V(M_1)\cup \cdots\cup V(M_p)\subseteq V(R_1')\cup \cdots\cup V(R_p').$$ 
Moreover, the sum of the lengths of $M_1\ll M_p$ is less than that of $R_1'\ll R_p'$,
since $u^1_1, v^1_1$ are not consecutive vertices of $R_1$. Hence if we take the linkage $L$ and replace each subpath $R_a'$ by $M_a$
for $1\le a\le p$, we obtain another linkage for the same problem instance using fewer vertices, contradicting that $L$ is minimum.
Thus the assumption immediately before (1) must have been false. This proves \ref{acceptable}.~\bbox

We deduce:

\begin{thm}\label{vertexorder}
Let $k,c,z$, $(G,s_1,t_1\ll s_k,t_k)$, $(C_1\ll C_c)$, and $L=(P_1\ll P_k)$ be as in \ref{acceptable}.
Then there is an enumeration $(v_1\ll v_n)$ of $V(L)\setminus \{s_1\ll s_k,t_1\ll t_k\}$, such that for $0\le h\le n$, 
if $B$ denotes $\{s_1\ll s_k\}\cup \{v_i:1\le i\le h\}$
and $A=V(L)\setminus B$, then 
\begin{itemize}
\item for $1\le a\le c$, $B\cap C_a$ is $C_a$-acceptable;
\item the $(B,A)$-wiggle number of each member of $L$ is at most $c(c-1)(z+1)+1$.
\end{itemize}
\end{thm}
\Proof 
Since $\{s_1\ll s_k\}$ is acceptable, repeated application of \ref{acceptable} implies that there is an 
enumeration $(v_1\ll v_n)$ of $V(L)\setminus \{s_1\ll s_k,t_1\ll t_k\}$, 
such that for $0\le h\le n$, if $B$ denotes $\{s_1\ll s_k\}\cup \{v_i:1\le i\le h\}$
and $A=V(L)\setminus B$, then
\begin{itemize}
\item for $1\le a\le c$, $B\cap C_a$ is $C_a$-acceptable;
\item for all distinct $a,b$ with $1\le a,b\le c$, and for $1\le i\le k$, the $(B\cap C_b, A\cap C_a)$-wiggle number of $P_i$ is
at most $z$.
\end{itemize}
We claim that this enumeration satisfies the theorem. For let $h,B,A$ be as in the theorem, let $1\le i\le k$, and let $t$
be the $(B,A)$-wiggle number of $P_i$. Consequently there are $t-1$ edges of $P_i$, say $a_1b_1\ll a_{t-1}b_{t-1}$,
such that $a_j\in A$ and $b_j\in B$ for $1\le j\le t-1$. For each such $j$, let $p_j,q_j\in \{1\ll c\}$ such that $a_j\in C_{p_j}$
and $b_j\in C_{q_j}$. Since $a_j\in A$ and $b_j\in B$, and $a_jb_j$ is a directed edge of $P_i$, it follows 
(since $B$ is $C_{p_j}$-acceptable) that $p_j\ne q_j$. There are only $c(c-1)$ possibilities for the pair $(p_j,q_j)$, 
and for each one of them, say $(p,q)$,
there are at most $z+1$ values of $j$ with $(p_j,q_j) = (p,q)$, since the $(B\cap C_q,A\cap C_p)$-wiggle number of $P_i$ is
at most $z$. Hence there are at most $c(c-1)(z+1)$ values of $j$ in total, and so $t\le c(c-1)(z+1)+1$, and this enumeration 
satisfies the theorem. This proves \ref{vertexorder}.~\bbox

\section{Enlarging on history}

In this section we define the sets $D_h$ and $\mathcal{D}$ discussed at the end of section 2, and use \ref{vertexorder}
to prove they have the properties we need.
Let $G$ be a digraph, and $(C_1\ll C_c)$ a clique-partition of $G$, and let $s$ be some positive integer. If $X$ is a sequence of vertices of $G$ we define
$V(X)$ to be the set of terms of $X$. Let $\mathcal{A}$ be a set of sequences of vertices of $G$.
We define $\mathcal{A}^+$ to be the set of $v\in V(G)$ such that for some $X\in \mathcal{A}$,
there exists $a\in \{1\ll c\}$ such that $\{v\}\cup V(X)\subseteq C_a$ and either
\begin{itemize}
\item $v\in V(X)$, or
\item $v\notin V(X)$ and $X$ has length $s$ and $v$ is adjacent from the last $s-1$ vertices of $X$.
\end{itemize}
(Thus, the order of the terms in $X\in \mathcal{A}$ does not matter, except it matters which term is first.)
Similarly, we define  $\mathcal{A}^-$ to be the set of vertices $v$ such that for some $X\in \mathcal{A}$,
there exists $a\in \{1\ll c\}$ such that $\{v\}\cup V(X)\subseteq C_a$ and either
\begin{itemize}
\item $v\in V(X)$, or
\item $v\notin V(X)$ and $X$ has length $s$ and $v$ is adjacent to the first $s-1$ vertices of $X$.
\end{itemize}

Now let $r,s,t\ge 0$ be integers. A subset $D$ of $V(G)$
is said to be {\em $(r,s,t)$-restricted} if there are sets $\mathcal{A}, \mathcal{B}$ of sequences of vertices of $G$, 
satisfying the following:
\begin{itemize}
\item every member of $\mathcal{A}$ and every member of $\mathcal{B}$ has length at most $s$;
\item $|\mathcal{A}|, |\mathcal{B}|\le r$; 
\item $|\mathcal{B}^+ \cap \mathcal{A}^-|\le t$; and
\item $\mathcal{B}^+ \setminus \mathcal{A}^-\subseteq D$, and $D\subseteq \mathcal{B}^+$.
\end{itemize}
Thus, for any constants $r,s,t$ there are only polynomially many $(r,s,t)$-restricted subsets $D$ of $V(G)$. For suitable
$r,s,t$ the set of all $(r,s,t)$-restricted subsets will be the set $\mathcal{D}$ that we need.

We observe:

\begin{thm}\label{shift}
Let $L$ be a minimum linkage for $(G,s_1,t_1\ll s_k,t_k)$, let $\ell\ge 3$, let $C$ be a clique, let $Q'$ be a subpath of some member of $L$, let $Q$ be a subpath of $Q'$, 
with $|C\cap V(Q)|\ge \ell$, and let $v\in C\setminus V(L)$ be adjacent from the last $\ell$ vertices of $Q$ in $C$. Then $v$ is adjacent from the
last $\ell$ vertices of $Q'$ in $C$.
\end{thm}
\Proof Let the vertices of $Q'$ in $C$
in order be $y_1\ll y_m$ say, and let the last $\ell$ vertices of $Q$ in $C$ be $x_1\ll x_{\ell}$ in order. 
Thus $m\ge \ell$,  since $x_1\ll x_{\ell}$ is a subsequence of $y_1\ll y_m$. Let $j\in \{m-\ell+1\ll m\}$.
We claim that $y_j$ is adjacent to $v$. For suppose not; then $y_j$ is different from all of $x_1\ll x_{\ell}$, and since $x_1\ll x_{\ell}$
are $\ell$ consecutive terms of the sequence $y_1\ll y_m$, and there are at most ${\ell}-1$ terms of this sequence after $y_j$,
it follows that $x_1\ll x_{\ell}$ all come before $y_j$. In particular, $x_1$ equals some $y_g$ where $g\le j-\ell$. Now $v$ is adjacent from
$x_1 = y_g$, and not adjacent from $y_j$. Since $v, y_j\in C$, it follows that $v$ is adjacent to $y_j$; but then replacing
the subpath of $P_i$ between $y_g,y_j$ by the path with three vertices $y_g\d v\d y_j$ contradicts that $L$ is a minimum linkage, since
$\ell\ge 3$. This proves \ref{shift}.~\bbox

If $P$ is a path of $G$, and $C\subseteq V(G)$ and $s$ is an integer, by the {\em first up-to-$s$ vertices of $P$ in $C$}
we mean the sequence which consists of the first $s$ vertices of $P$ in $C$, if there are $s$ such vertices, and otherwise 
the sequence consisting of all vertices of $P$ in $C$, in either case in their order in $P$. We define the {\em last up-to-$s$} similarly.

Next we define the sets $D_h$.
Let $(G,s_1,t_1\ll s_k,t_k)$ be a problem instance, where
$G$ admits a clique-partition $(C_1\ll C_c)$, and let $L$ be a 
linkage for this problem instance. Let $s = k^2+k+3$.
Let $(v_1\ll v_n)$ be an enumeration of $V(L)\setminus \{s_1\ll s_k,t_1\ll t_k\}$, and for $0\le h\le n$, let $B_h$ denote $\{s_1\ll s_k\}\cup \{v_i:1\le i\le h\}$
and $A_h=V(L)\setminus B_h$. For $0\le h\le n$, let $J_h$ be the set of edges of $G$ that belong to a member of $L$
and have one end in $A_h$ and the other in $B_h$.
The union of the members of $L$ is a digraph consisting of $k$
disjoint paths, and if we delete $J_h$ from this digraph, we obtain a digraph which is 
also a disjoint union of paths, each with vertex set included in one of $A_h,B_h$. Let $\mathcal{Q}_h$ 
be the set of these paths which are included in $B_h$, and $\mathcal{R}_h$ the set included in $A_h$.
Let $\mathcal{A}_h$
be the set of all sequences $X$ such that for some $R\in \mathcal{R}_h$ and some $a\in \{1\ll c\}$, $X$ is the first up-to-$s$ vertices of $R$ in $C_a$.
Similarly, let $\mathcal{B}_h$ be the set of all sequences
$X$ such that for some $Q\in \mathcal{Q}_h$ and $a\in \{1\ll c\}$, $X$ is the last up-to-$s$ vertices of $Q$ in $C_a$.
We define $D_h=(\mathcal{B}_h^+\setminus \mathcal{A}_h^-)\cup B_h$.

We claim:

\begin{thm}\label{separation}
Let $(G,s_1,t_1\ll s_k,t_k)$,
$(C_1\ll C_c)$, $L$,
$(v_1\ll v_n)$, and $A_h, B_h\;(0\le h\le n)$ be as above, and let $w\ge 0$.
Suppose that 
\begin{itemize}
\item $L$ is a minimum linkage for $(G,s_1,t_1\ll s_k,t_k)$;
\item for $1\le a\le c$, $B_h\cap C_a$ is $C_a$-acceptable; and
\item the $(A_h,B_h)$-wiggle number of each member of $L$ is at most $w$.
\end{itemize}
Let $r=ckw$, $s = k^2+k+3$, $t = 2cskw+c(2w+1)k^2(k^2+k+1)$, and for $0\le h\le n$ let $D_h$
be as above. 
Then 
\begin{itemize}
\item[{\rm(a)}] $B_h\subseteq D_h$ and $A_h\cap D_h = \emptyset$ for $0\le h\le n$;
\item[{\rm(b)}] $D_h\subseteq D_{h+1}$ for $0\le h<n$; and
\item[{\rm(c)}] $D_h$ is $(r,s,t)$-restricted for $0\le h\le n$.
\end{itemize}
\end{thm}
\Proof
Let $0\le h\le n$.
Since the $(A_h,B_h)$-wiggle number of each member of $L$ is at most $w$, there are at most $2w-1$
edges of each member of $L$ in $J_h$, and the sets $\mathcal{Q}_h, \mathcal{R}_h$ defined in the 
definition of $D_h$ both have at most $kw$ members.
Thus the sets $\mathcal{A}_h, \mathcal{B}_h$
both have cardinality at most $ckw=r$.  
\\
\\
(1) {\em $|\mathcal{B}^+_h \cap \mathcal{A}^-_h|\le t$.}
\\
\\
There are at most $kw$ choices for $Q\in \mathcal{Q}_h$, and for each there are at most $c$ choices for the sequence of the
last up-to-$s$ vertices of $C_a$ in $Q$, one for each $a\in \{1\ll c\}$; and each such sequence has at most $s$ terms. Thus there
are at most $cskw$ vertices which belong to the sequence of the last up-to-$s$ members in some $C_a$ of some path $Q\in \mathcal{Q}_h$.
Similarly there are at most $cskw$ vertices that belong to the first up-to-$s$ members of some $C_a$ in some $R\in \mathcal{R}_h$,
a total of at most $2cskw$.
For every other vertex $v\in \mathcal{B}^+_h \cap \mathcal{A}^-_h$, choose 
$a\in \{1\ll c\}$ such that $v\in C_a$; then
\begin{itemize}
\item[($\ast$)] there exists $Q\in \mathcal{Q}_h$ such that 
$|C_a\cap V(Q)|\ge s$, and $v$ is not among the last $s-1$ vertices of $C_a$ in $Q$, and is adjacent from each of the last $s-1$; and 
there exists $R\in \mathcal{R}_h$ such that 
$|C_a\cap V(R)|\ge s$, and $v$ is not among the first $s-1$ vertices of $C_a$ in $R$, and $v$ is adjacent to each of
the first $s-1$.
\end{itemize}
Let us fix $a\in \{1\ll c\}$, and count the number of vertices $v\in C_a$ satisfying ($\ast$).
Such vertices $v$ might belong to $A_h$ or to $B_h$ or to neither, and we count
the three types separately. First, suppose that $v\in A_h$. Thus $v$ belongs to some member $P_i$ of $L$; and there are only $k$
choices for $i$. For each $Q\in \mathcal{Q}_h$ containing at least $s$ vertices in $C_a$, let $X$ be the set of the last $s-1$ 
such vertices of $C_a$ in $Q$; there are at most $k^2+k+1$ vertices in $C_a\cap V(P_i)\cap A_h$ that are adjacent from every
vertex in $X$, 
since $B_h$ is $C_a$-acceptable.
Since there are only $kw$ choices of $Q$, it follows that there are at most $wk(k^2+k+1)$ vertices $v\in C_a\cap V(P_i)\cap A_h$
satisfying ($\ast$); and summing over $1\le i\le k$, we deduce there are
at most $wk^2(k^2+k+1)$ vertices $v\in C_a\cap A_h$
satisfying ($\ast$). Similarly there are at most that many in $C_a\cap B_h$.

Finally, we must count the number of $v\in C_a\setminus V(L)$ satisfy ($\ast$).
By \ref{shift}, if $v\in C_a\setminus V(L)$ and is adjacent from the last $s-1$ vertices of $C_a\cap B_h$ in some subpath of $P_i$, then $v$ is also adjacent from the
last $s-1$ vertices of $C_a\cap B_h$ in $P_i$. We deduce that if $v\in C_a\setminus V(L)$, and $v$ satisfies 
($\ast$), then for some $i,j\in \{1\ll k\}$, $v$ is adjacent from the last $s-1$ vertices of $P_i$ in
$B_h\cap C_a$, and adjacent to the first $s-1$ vertices of $P_j$ in $A_h\cap C_a$ (similarly). For any choice of $i,j$
there are at most $k^2+k+1$ such vertices $v$, because $B_h$ is $C_a$-acceptable. (This is where we use paths of length two in 
the definition of a planar matching.) Consequently there are at most $k^2(k^2+k+1)$ such vertices $v\in C_a\setminus V(L)$ total.

Altogether, we have shown that there are at most $2wk^2(k^2+k+1) + k^2(k^2+k+1)$ vertices $v\in C_a$
satisfying ($\ast$), and summing over $a\in \{1\ll c\}$ and adding back the $2cskw$ from the start of the argument, 
the claim follows. This proves (1).
\\
\\
(2) {\em $A_h\subseteq \mathcal{A}^-_h$ and $B_h\subseteq \mathcal{B}^+_h$.}
\\
\\
Let $v\in A_h$; then $v$ belongs to $V(L)$, and hence to some member of $L$, and therefore to some member of $\mathcal{R}_h$, say $R$.
Choose $a\in \{1\ll c\}$ with $v\in C_a$, and let $X\in \mathcal{A}_h$ be the sequence of the first up-to-$s$ vertices of $R$
in $C_a$.
If $v\in V(X)$ then $v\in \mathcal{A}^-_h$ as required, so we may assume not; and so there are more than $s$
vertices of $R$ in $C_a$, and $X$ has exactly $s$ terms. Let the vertices of $X$ be $x_1\ll x_s$ in order. Then $x_1\ll x_s,v\in C_a$,
and $x_i$ is not adjacent to $v$ for $1\le i<s$, since $R$ is a minimal path of $G$ (because the members of $L$ are minimal paths).
Thus $v$ is adjacent to $x_1\ll x_{s-1}$, and hence $v\in \mathcal{A}^-_h$ as required. Similarly $B_h\subseteq \mathcal{B}^+_h$. This proves
(2).
\\
\\
(3) {\em $B_h\subseteq D_h$, $A_h\cap D_h = \emptyset$ and $D_h\subseteq \mathcal{B}^+_h$.}
\\
\\
We recall that $D_h = (\mathcal{B}^+_h\setminus \mathcal{A}^-_h)\cup B_h$. Consequently $B_h\subseteq D_h$, and from (2),
$D_h\subseteq \mathcal{B}^+_h$. Since $A_h\cap B_h=\emptyset$ and $A_h\subseteq \mathcal{A}^-_h$ by (2), it follows that 
$A_h\cap D_h = \emptyset$. This proves (3).

\bigskip

Assertion (a) of the theorem follows from (3), and (c) from (1) and (3). We still
need to show (b). Let $0\le h<n$; we must show that $D_h\subseteq D_{h+1}$. Let $v\in D_h$; we will
show that $v\in D_{h+1}$. If $v\in B_h$
then $v\in B_{h+1}\subseteq D_{h+1}$ as required, so we assume that $v\notin B_h$. Certainly $v\notin A_h$ by (3)
since $v\in D_h$; so $v\notin V(L)$. Since $v\in D_h$, it follows that $v\in \mathcal{B}_h^+\setminus \mathcal{A}_h^-$,
and in particular there exist $Q\in \mathcal{Q}_h$ and $a\in \{1\ll c\}$ such that $v\in C_a$, and $|C_a\cap V(Q)|\ge s$,
and $v$ is adjacent from the last $s-1$ vertices of $C_a$ in $Q$.
Let $Q$ be a subpath of $P_i\in L$, and let $Q'$ 
be the maximal subpath of $P_i$ including $Q$ such that all its vertices are in $B_{h+1}$. By \ref{shift}, $v$ is adjacent from the
last $s-1$ vertices of $C_a$ in $Q'$; and so $v\in \mathcal{B}_{h+1}^+$. It remains
to show that $v\notin \mathcal{A}_{h+1}^-$; so, suppose it is. Then by the same argument with $A_h$ exchanged with $B_{h+1}$,
and $A_{h+1}$ exchanged with $B_h$ (and $h,h+1$ exchanged) it follows that $v\in \mathcal{A}_{h}^-$, a contradiction. This
proves assertion (b) of the theorem, and so completes the proof of \ref{separation}.~\bbox

\section{The auxiliary digraph}

Let $k,c\ge 1$ and let $r,s,t, w$ be as in \ref{separation}.  
Let $(G,s_1,t_1\ll s_k, t_k)$ be a problem instance, where $G$ admits a clique-partition
$(C_1\ll C_c)$. Let $\mathcal{D}$ be the set of all $(r,s,t)$-restricted subsets of $V(G)$. A {\em coloured edge} means
a pair $(e,i)$ where $e\in E(G)$ and $1\le i\le k$, and we will abuse this terminology a little, speaking of coloured edges as though they are edges
(for instance, we speak of the {\em head} of a coloured edge $(e,i)$ meaning the head of $e$,
and so on). We call $i$ the {\em colour} of the coloured edge. Let $\mathcal{E}$ be the set of all 
sets $Y$ of coloured edges of cardinality at most $2w-1$, such that 
\begin{itemize}
\item no two members of $Y$ have the same head or 
the same tail, and 
\item every two members of $Y$ that share an end have the same colour; 
\item no coloured edge in $Y$ has head in $\{s_1\ll s_k\}$ or tail in $\{t_1\ll t_k\}$; and 
\item for $1\le i\le k$, every coloured edge with tail $s_i$ has colour $i$, and every coloured edge with head $t_i$ has colour $i$.
\end{itemize}
The auxiliary digraph $H$ will have vertex set all pairs
$(Y,D)$ where $Y\in \mathcal{E}$ and $D\in \mathcal{D}$ and every coloured edge in $Y$ has exactly one end in $D$.

Now we define its adjacency. Let $(Y,D), (Y', D')\in V(H)$ be distinct. 
We say that $(Y,D)$ is adjacent to $(Y', D')$ in $H$ if $D\subseteq D'$, and there are exactly two coloured edges that belong to 
$(Y\setminus Y')\cup (Y'\setminus Y)$, and they form a two-edge path with middle vertex in $D'\setminus D$.

That defines $H$. Now let $S_0$ be the set of all vertices $(Y,D)$ of $H$ such that $|Y|=k$ and every coloured edge in $Y$
has tail in $\{s_1\ll s_k\}$, and let $T_0$ be the set of all $(Y,D)$ such that $|Y|=k$ and every coloured edge in $Y$ has
head in $\{t_1\ll t_k\}$. We claim:

\begin{thm}\label{equivalence}
Let $k,c\ge 1$, and let $r,s,t,w$ be as in \ref{separation}. Let $(G,s_1,t_1\ll s_k, t_k)$ be a problem instance, where $G$ admits a clique-partition
$(C_1\ll C_c)$, and let $\mathcal{D},\mathcal{E}$ and $H,S_0,T_0$ be as above. Then there is a path in $H$ from a vertex in $S_0$ to 
a vertex in $T_0$ if and only if there is a linkage for $(G,s_1,t_1\ll s_k, t_k)$.
\end{thm}
\Proof
We observe first (the proof is clear and we omit it):
\\
\\
(1) {\em If there is a directed path in $H$ from $(Y,D)$ to $(Y',D')$ then $D\subseteq D'$.}

\bigskip

Suppose that there is a path in $H$ from $S_0$ to $T_0$, with vertices
$(Y_1,D_1)\ll (Y_n,D_n)$ say, in order. 
Let $Y = Y_1\cup\cdots\cup Y_n$; thus $Y$ is a set of coloured edges. We need to show that
$Y$ includes the edge set of a linkage for $(G,s_1,t_1\ll s_k, t_k)$. 
\\
\\
(2) {\em If $(e,i)\in Y$ and $1\le h\le n$, and exactly one end of $e$ is in $D_h$, then $(e,i)\in Y_h$.}
\\
\\
For let $(e,i)\in Y_{h'}$ where $1\le h'\le n$; and choose $h'$ with $|h-h'|$ minimum. 
Let $e$ have ends $u,v$, and suppose that $h\ne h'$. Now exactly one end of $e$ belongs to $D_{h'}$, and exactly one to $D_h$, and one of $D_h, D_{h'}$
is a subset of the other by (1); and so it is the same
end $v$ say of $e$ that lies in both $D_h, D_{h'}$, and $u$ belongs to neither of them. 
If $h'<h$ let $h'' = h'+1$, and if $h'>h$ let $h'' = h'-1$. Since one of $D_h,D_{h'}$ is a subset of $D_{h''}$ (by (1)) 
it follows that
$v\in D_{h'}$; and since $D_{h''}$ is a subset of one of $D_h,D_{h'}$ (by (1) again) it follows that $u\notin D_{h''}$. But this contradicts
the minimality of $|h'-h|$. Consequently $h=h'$, and the claim holds. This proves (2).
\\
\\
(3) {\em If $(e,i), (e',i')\in Y$ share an end then $i=i'$, and these edges form a two-edge path.}
\\
\\
Choose $h$ with $1\le h\le n$ such that $(e,i)\in Y_{h}$, and choose $h'$ similarly for $(e',i')$; and
in addition choose $h,h'$ with $|h-h'|$ minimum. If $h= h'$, then $(e,i), (e',i')$ are both in $Y_{h}$,
and since they share an end, it follows that $i=i'$ and the edges form a two-edge path, from the definition of $\mathcal{E}$.
Thus we may assume that $h<h'$. Now $D_h\subseteq D_{h'}$ by (1), and since one end of $e$ belongs to $D_h$, it follows that at least 
one end of $e$ is in $D_{h'}$; and since $(e,i)\notin Y_{h'}$, (2) implies that both ends of $e$ are in $D_{h'}$. 
Similarly, neither end of $e'$ is in $D_h$; and so the common end of $e,e'$ belongs to $D_{h'}\setminus D_h$. Let $e$ have
ends $u,v$, and let $e'$ have ends $v, w$, where $u\in D_h$, $v\in D_{h'}\setminus D_h$, and $w\notin D_{h'}$.
Now $u\in D_{h+1}$ since $D_h\subseteq D_{h+1}$; and since $(e,i)\notin Y_{h+1}$, (2) implies that $v\in D_{h+1}$. Consequently
$(e',i')\in Y_{i+1}$ by (2), and hence $h' = h+1$ from the minimality of $|h-h'|$. But now the claim follows from the definition of
adjacency in $H$. This proves (3).
\\
\\
(4) {\em Every vertex of $G$ incident with exactly one coloured edge in $Y$ belongs to $\{s_1 \ll s_k, t_1\ll t_k\}$.}
\\
\\
For suppose that $v\in V(G)\setminus \{s_1 \ll s_k, t_1\ll t_k\}$ is incident with exactly one coloured edge $(e,i)\in Y$.
Let the other end of $e$ be $u$. There are three cases, depending whether $u\in \{s_1\ll s_k\}$, $u\in \{t_1\ll t_k\}$, or neither.
Suppose first that $u\in \{s_1\ll s_k\}$. Then $(e,i)\in Y_1$; and $(e,i)\notin Y_n$ since $v\notin \{t_1\ll t_k\}$.
Choose $h<n$ maximum such that $(e,i)\in Y_h$. From the maximality of $h$, $(e,i)\notin Y_{h+1}$, and so by (2)
$u,v\in D_{h+1}$. From the definition of adjacency in $H$, there is another edge $(f,i)\in Y$ forming a two-edge path with
$(e,i)$, such that the common end of $e,f$ is not in $D_h$. But $u\in \{s_1\ll s_k\}\subseteq D_h$, and $v$ is not incident
with any other edge in $Y$, a contradiction. The argument is analogous if $u\in \{t_1\ll t_k\}$ and we omit it.
Finally suppose that $u\notin \{s_1 \ll s_k, t_1\ll t_k\}$. Thus $(e,i)\notin Y_1, Y_n$; choose $h<h'<h''$ with $h''-h$ minimum such
that $(e,i)\notin Y_h\cup Y_{h''}$ and $(e,i)\in Y_{h'}$. From the minimality of $h''-h$ it follows that
$(e,i) \in Y_{h+1}$; and from the definition of adjacency in $H$, there is an edge $(f,i)$ of $Y$ that makes a two-edge
path with $(e,i)$, such that the common end of $e,f$ is in $D_{h+1}\setminus D_h$. Since $v$ is not incident
with any other edge in $Y$, this common end is $u$, so $u\in D_{h+1}\setminus D_h$. But similarly, $u\in D_{h''}\setminus D_{h''-1}$, and this is
impossible since $D_{h+1}\subseteq D_{h''-1}$ by (1). This proves (4).

\bigskip

From (3), no three edges in $Y$ share an end (because this end would be the head of two of them or the tail of two,
contrary to (3)). Thus the digraph formed by $Y$ is the disjoint union of directed paths and directed cycles, and we call
these ``components'' of $Y$. The edges in each component all have the same colour, by (3). 
Each path component has first vertex in $\{s_1\ll s_k\}$ and last vertex in $\{t_1\ll t_k\}$,
by (4). Moreover,
for $1\le i\le k$, some edge in $Y_1\subseteq Y$ has tail $s_i$ and colour $i$ (from the definition of $S_0$); and 
since no edge in $Y$ has head $s_i$,
it follows that $s_i$ is the first vertex of a path component $P_i$ of $Y$ in which all edges
have colour $i$. The last vertex of this component is in $\{t_1\ll t_k\}$, and is therefore $t_i$ since the last edge has colour $i$.
Consequently $(P_1\ll P_k)$ is a linkage for $(G,s_1,t_1\ll s_k, t_k)$. This proves the ``only if'' part of the theorem.

Now we turn to the ``if'' part. We assume there is a linkage for $(G,s_1,t_1\ll s_k, t_k)$, and must prove there is a path in $G$
from $S_0$ to $T_0$. Let $L = (P_1\ll P_k)$ be a minimum linkage. Let $v_1\ll v_n$ be as in \ref{vertexorder}.
For $0\le h\le n$, let $B_h = \{s_1\ll s_k\}\cup \{v_i:1\le i\le h\}$ and $A_h = V(L)\setminus B_h$; and let $D_h$ be as 
defined immediately before \ref{separation}. Let $J_h$ be the set of edges of $P_1\cup\cdots\cup P_k$
with one end in $A_h$ and the other in $B_h$, and let $Y_h=\{(e,i):e\in J_h\cap E(P_i)\}$. We claim that
\begin{itemize}
\item $(Y_h,D_h)\in V(H)$ for $0\le h\le n$;
\item for $0\le h<n$, $(Y_h, D_h)$ is adjacent in $H$ to $(Y_{h+1}, D_{h+1})$; and 
\item $(Y_0,D_0)\in S_0$, and $(Y_n, D_n)\in T_0$.
\end{itemize}
To see the first claim, note that $D_h$ is $(r,s,t)$-restricted by \ref{separation}; and $Y_h\in \mathcal{E}$ since $L$ is a linkage.
Also the third claim follows. For the second, let $0\le h<n$. By \ref{separation}, $D_h\subseteq D_{h+1}$. Let $(e,i)$
be a coloured edge that belongs to exactly one of $Y_h, Y_{h+1}$. 
It follows that $e\in E(P_i)$, and hence has both ends in 
$V(L)$; and since $e$ belongs to exactly one of $J_h, J_{h+1}$, some end $v$ of $e$ belongs to $D_{h+1}\setminus D_h$. 
Thus $v\in D_{h+1}\cap V(L)=B_{h+1}$, and $v\in V(L)\setminus D_h = A_h$, by \ref{separation}. Hence $v = v_{h+1}$
since $B_{h+1}=B_h\cup \{v_{h+1}\}$. Now $v\notin \{s_1\ll s_k,t_1\ll t_k\}$, and so there is a two-edge subpath $Q$ of $P_i$
with middle vertex $v$. Since $B_{h+1}=B_h\cup \{v\}$, it follows that the other edge of $Q$ also belongs to 
exactly one of $J_h, J_{h+1}$; and no other edges have this property, since we have shown that every edge in exactly one of 
$J_h, J_{h+1}$ is incident with $v=v_{h+1}$, and no other edges in $P_1\cup\cdots\cup P_k$ are incident with $v$. This completes
the proof of the second bullet above, and so proves the ``if'' half of the theorem, and hence completes the proof of \ref{equivalence}.~\bbox

Let us figure out the running time. Checking whether the path in $H$ exists can be done in time $O(|V(H)|^2)$ 
(for instance by breadth-first search), which is
also the time needed to construct $H$; so we just need to estimate $|V(H)|$. We recall that
$z=c(c(k^2+k+1)+k+2)$, $w= c(c-1)(z+1)+1$, $r=ckw$, $s = k^2+k+3$, 
and $t = 2cskw+c(2w+1)k^2(k^2+k+1)$; and let $n=|V(G)|$.
Now $H$ has at most $|\mathcal{D}|\cdot |\mathcal{E}|$ vertices, and  
$|\mathcal{D}|\le n^{2rs}2^t$, and $|\mathcal{E}|\le (n^2k)^{2w-1}$. Hence $|V(H)|^2= O(n^{4rs+8w})$,
and this exponent is about $4(ck)^{5}$ for large $c,k$.

Finally, we remark that
every $p$-vertex path from $S_0$ to $T_0$ in $H$ gives a linkage in $G$ using at most $p-1+2k$ vertices; and 
every minimum linkage in $G$ with $p-1+2k$ vertices gives a $p$-vertex path in $H$. Thus if we check for the shortest path 
in $H$ from $S_0$ to $T_0$, we can determine the minimum number of vertices in a linkage for $(G,s_1,t_1\ll s_k,t_k)$, as mentioned
in the introduction.

\end{document}